\documentclass[12pt,twoside]{article}
\usepackage[centertags]{amsmath}
\usepackage{amsfonts}
\usepackage{amssymb}
\usepackage{amsthm}
\usepackage{newlfont}
\usepackage{color}

\date{\color{green}  2013 November}
\def \n {\noindent}
\usepackage{fancyhdr}
\newcommand{\hf}{\hfill $\diamondsuit$}
\begin{document}
\centerline{\bf \bf{\color{green} :::::::}}

\begin{center}
{\bf{\color{red}{\it January  2014}}}\\

***\\

{\Large{\color{blue}{\bf On the complete indeterminacy and the chaoticity of generalized Heun's operator in Bargmann space}}}\\
\end{center}

\begin{center}
{\bf Abdelkader INTISSAR  $^{(*)}$ $^{(**)}$}\\
\end{center}

\begin{center}
\scriptsize{(*)Equiped'Analyse spectrale, UMR-CNRS n: 6134, Universit\'e de
Corse, Quartier Grossetti, 20 250 Cort\'e-France \\
 T\'el: 00 33 (0) 4 95 45 00 33\\
 Fax: 00 33 (0) 4 95 45 00 33\\
 e.mail:intissar@univ-corse.fr\\
\parskip = 3pt
(**)Le Prador,129 rue du commandant Rolland, 13008 Marseille-France}
\end{center}

\quad\\

\begin{center}

\fbox{\rule[-0.4cm]{0cm}{1cm} \n {\color{red}{\bf Abstract}}}\\

\end{center}

{\it \small{In Communications in Mathematical Physics, no. 199, (1998), we have considered the Heun operator  $\displaystyle{ H = a^* (a + a^*)a}$  acting on Bargmann space where $a$ and $a^{*}$ are the standard Bose annihilation and creation operators satisfying the commutation relation $[a, a^{*}] = I$.\\

We have used the boundary conditions at infinity to give a description of all maximal
dissipative extensions in Bargmann space of the minimal Heun's operator $H$. The characteristic functions of the dissipative extensions  have be computed and some completeness theorems have be obtained for the system of generalized eigenvectors of this operator. \\

\quad In this paper we study the deficiency numbers of the generalized Heun's operator
 $\displaystyle{ H^{p,m} = a^{*^{p}} (a^{m} + a^{*^{m}})a^{p}; (p, m=1, 2, .....)}$
acting on Bargmann space. In particular, here we find some conditions on the parameters $p$ and $m$ for that $\displaystyle{ H^{p,m}}$ to be completely indeterminate. It follows from these conditions that $\displaystyle{ H^{p,m}}$ is entire of the type minimal. And we show that $\displaystyle{ H^{p,m}}$ and $\displaystyle{ H^{p,m}+ H^{*^{p,m}}}$ (where $H^{*^{p,m}}$ is the adjoint of the $H^{p,m}$) are connected at the chaotic operators. We will give a description of all maximal dissipative extensions and all selfadjoint extensions of the minimal generalized Heun's operator $H^{p,m}$ acting on Bargmann space in separate paper.\\

 }}\\

 \begin{center}
\Large{\bf{\color{blue} 1. Introduction and preliminaries results}}
\end{center}

 Let  $\mathbb{B}$ be the Bargmann space {\bf[3]} defined as a subspace of the space $O( \mathbb{C})$ of holomorphic functions on $\mathbb{C}$, given by\\

\n $\mathbb{B} = \{\phi\in O(\mathbb{C}) ; < \phi, \phi > < \infty \} $ $\hfill {(1.1) } $\\

\n where the paring \\

\n $<\phi,\psi> = \displaystyle{\int_{\mathbb{C}}}\displaystyle{\phi(z)\overline{\psi(z)}e^{-\mid z\mid^{2}}dxdy}$ $\forall$ $\phi,\psi \in O(\mathbb{C})$ $\hfill {(1.2) } $\\

\n and $dxdy$ is Lebesgue measure on $\mathbb{C}$.\\

\n This space  with $\mid\mid \phi \mid\mid = \sqrt{<\phi, \phi>}$ is a Hilbert space and\\ $\displaystyle{e_{k}(z) = \frac{z^k}{\sqrt{k!}}; k = 0, 1, ....}$ is complete  orthonormal basis of $\mathbb{B}$.\\

In this representation, the standard Bose annihilation and creation operators are defined by\\

$\left\{\begin{array}[c]{l}a\phi(z) = \quad \phi^{'}(z)\\ with \quad maximal \quad domain\\

D(a) = \{\phi \in \mathbb{B} \quad such \quad that \quad  a\phi \in \mathbb{B}\} \\ \end{array}\right.
\hfill { } (1.3) $\\
\quad\\

$\left\{\begin{array}[c]{l}a^{*}\phi(z) = \quad z\phi(z)\\ with \quad maximal \quad domain\\

D(a^{*}) = \{\phi \in \mathbb{B} \quad such \quad that \quad  a^{*}\phi \in \mathbb{B}\} \\ \end{array}\right.\hfill { }  (1.4)$\\
\quad\\

Accordingly, for the operator  $ H^{p,m} = a^{*^{p}} (a^{m} + a^{*^{m}})a^{p}; (p, m =1, 2, ..... ) $ we have\\

$\left\{\begin{array}[c]{l}H^{p,m} \phi(z)=   z^{p}\phi^{(p+m)}(z) + z^{p+m}\phi^{(p)}(z) \quad \quad \quad  \quad \quad \quad\\ with \quad maximal \quad domain\\

D(H_{max}^{p,m}) = \{\phi \in \mathbb{B} ; H^{p,m}\phi \in \mathbb{B}\} \\ \end{array}\right.\hfill { } (1.5)$.\\

It follows from (1.3) and (1.4) that the action of the operator $a^{*^{r}}a^{s}$\\ ($r \in \mathbb{N}$, $s \in \mathbb{N}$) on an element $\phi \in \mathbb{B}$; $\displaystyle{\phi(z) = \sum_{k=0}^{\infty}a_{k}e_{k}(z)}$ is given by:\\

$\displaystyle{a^{*^{r}}a^{s}\phi(z) }$  = $\displaystyle{\sum_{k=0}^{\infty}k(k-1).....(k-s+1)a_{k}\frac{z^{k+r-s}}{\sqrt{k!}}}$\\.\quad \quad \quad \quad \quad =
$\displaystyle{\sum_{k=0}^{\infty}(k+s-r)(k+s-r-1).....(k-r+1)a_{k+s-r}\frac{z^{k}}{\sqrt{(k+s-r)!}}}$\\

If $s \geq r$, we have $\displaystyle{\sqrt{(k + s -r)!}= \sqrt{k!} \sqrt{(k + 1)!}..... \sqrt{(k + s -r)!}}$ this implies that\\

$\displaystyle{a^{*^{r}}a^{s}\phi(z) }$  = $\displaystyle{\sum_{k=0}^{\infty}\{\sqrt{(k + s -r)}.\sqrt{(k + s -r-1)}.....\sqrt{(k + 1)}.......(k-r+1) a_{k+s-r}\}e_{k}(z)}$\\

Let $\displaystyle{u_{k} = \sqrt{(k + s -r)}.\sqrt{(k + s -r-1)}.....\sqrt{(k + 1)}.......(k-r+1)}$ then we give the bellow obvious lemma that we will use in the following of this paper\\

{\bf Lemma 1.1}\\

i) if $s \geq r$ then $\displaystyle{u_{k}\sim k^{\frac{s+r}{2}}}$\\

ii) Also, if $s \leq r$ then $\displaystyle{u_{k}\sim k^{\frac{s+r}{2}}}$\\

iii) For $r = p$ and $s = p + m$ then $\displaystyle{u_{k}\sim k^{\frac{2p+m}{2}}}$\\

{\bf Proof}\\

  As $\displaystyle{u_{k}\sim k^{\frac{s-r}{2}}k^{r}}$ then $\displaystyle{u_{k}\sim k^{\frac{s+r}{2}}}$. In a similar manner, we obtain $u_{k}$ in ii) when $s \leq r$ and iii) is obvious particular case of i). The  proof of Lemma 1.1 is complete.{\color{blue}\hf}\\

 Now the action of the operator $H^{p,m}$ on an element $\phi \in \mathbb{B}$ is given by:\\

$\displaystyle{H^{p,m}\phi(z) = \sum_{k=0}^{\infty}k(k-1)........(k-p - m +1)a_{k}\frac{z^{k-m}}{\sqrt{k!}}\quad +}$ \\

$\displaystyle{\quad \quad \quad \quad \quad\sum_{k=0}^{\infty}k(k-1).......(k - p + 1)a_{k}\frac{z^{k+m}}{\sqrt{k!}}}$\\

As $H^{p,m}$ is polynomial in $(a^{*}, a)$ of degree $2p + m$, we give the next lemma for study in separate paper an perturbation of  $H^{p,m}$ by selfadjoint operators of the form $(a^{*^{j}}a^{j})$ with $j > p + \frac{m}{2}$, in particulary we will show the non chaoticity of this perturbation and we will give a complete spectral analysis for the genre of these operators.\\

{\bf Lemma 1.2}\\

For all $j > p + \frac{m}{2}$ the following statement holds:\\

$\displaystyle{\forall \quad \epsilon > 0\quad \exists\quad C_{\epsilon} > 0 }$ such that\\

 $ \displaystyle{ \forall \quad \phi \in D(a^{2j})}$,  $\displaystyle{\mid < H^{p,m}\phi , \phi > \mid\leq \epsilon \mid\mid a^{j}\phi\mid\mid^{2} + C_{\epsilon}\mid\mid \phi \mid\mid^{2}}$.\\

{\bf Proof}\\

Let $\displaystyle{\phi \in D(a^{2j})}$  and $\displaystyle{u_{k} = \sqrt{(k + m)}.\sqrt{(k + m -1)}.....\sqrt{(k + 1)}.......(k-p+1)}$ then\\

$\displaystyle{a^{*^{p}}a^{p+m}\phi, \phi > = \sum_{k=0}^{\infty}u_{k}a_{k+m}\bar{a}_{k}}$\\

$\displaystyle{\leq \sum_{k=0}^{\infty}u_{k}\mid a_{k+m}\mid \mid \bar{a}_{k}\mid }$\\

$\displaystyle{\leq \frac{1}{2}\sum_{k=0}^{\infty}u_{k}\mid a_{k+m}\mid^{2} + \frac{1}{2}\sum_{k=0}^{\infty}u_{k}\mid a_{k}\mid^{2} }$\\

$\displaystyle{\leq \frac{1}{2}\sum_{k=0}^{\infty}(u_{k-m} + u_{k})\mid a_{k}\mid^{2}}$\\

By virtue of iii) of lemma 1.1, we have $\displaystyle{u_{k}\sim k^{\frac{2p+m}{2}}}$ then there exist $c_{0} > 0$ and $c_{1} > 0$ such that
$\displaystyle{u_{k} \leq c_{0} + c_{1} k^{\frac{2p+m}{2}}}$\\

Now, for $j > p + \frac{m}{2}$ we apply the Young's inequality to get\\

$\displaystyle{\forall\quad \delta > 0\quad \exists\quad c_{\delta} > 0; k^{\frac{2p+m}{2}}\leq \delta k^{j} + c_{\delta}}$\\

this implies that\\

$\displaystyle{a^{*^{p}}a^{p+m}\phi, \phi > \leq c_{1}\delta \sum_{k=0}^{\infty}k^{j}\mid a_{k}\mid^{2} + (c_{\delta} + c_{0})\sum_{k=0}^{\infty}\mid a_{k}\mid^{2}}$\\

and\\

$\displaystyle{\forall \quad \epsilon > 0\quad \exists\quad C_{\epsilon} > 0 }$ such that
$\displaystyle{ \forall \quad \phi \in D(a^{2j})}$, \\

$\displaystyle{\mid < a^{*^{p}}a^{p+m}\phi , \phi > \mid\leq \epsilon \mid <a^{*^{j}}a^{j}\phi, \phi >\mid + C_{\epsilon}\mid\mid \phi \mid\mid^{2}}$.\\

As $\displaystyle{< a^{*^{p+m}}a^{p}\phi , \phi > = < \phi , a^{*^{p}}a^{p+m}\phi >}$ then we get\\

$\displaystyle{\forall \quad \epsilon > 0\quad \exists\quad C_{\epsilon} > 0 }$ such that
$\displaystyle{ \forall \quad \phi \in D(a^{2j})}$,\\

$\displaystyle{\mid < H^{p,m}\phi , \phi >\mid \leq \epsilon \mid\mid a^{j}\phi\mid\mid^{2} + C_{\epsilon}\mid\mid \phi \mid\mid^{2}}$. \\

The  proof of the Lemma 1.2 is complete. {\color{blue}\hf}\\

Now, the differential operations $a^{*}$
and $a$ act on the functions $e_{k}$ according to the formulas\\

$\displaystyle{a^{*}e_{k} = \sqrt{k+1}e_{k+1}, ae_{k} = \sqrt{k}e_{k-1}; e_{-1} = 0, k = 0, 1 , . . .}$ $\hfill { } (1.6)$\\

It follows from (1.6) that\\

$\left\{\begin{array}[c]{l}\displaystyle{H^{p,m}e_{k} = 0} \quad if \quad k < p\\\quad\\
\displaystyle{H^{p,m}e_{k} = \frac{\sqrt{k!(k + m)!}}{(k - p)!}e_{k+m}}\quad if \quad p\leq k < p+m  \end{array}\right.\hfill { } (1.7)$\\

and\\

$\displaystyle{H^{p,m}e_{k} = \frac{\sqrt{k!(k - m)!}}{(k - p - m)!}e_{k-m}}$ + $\displaystyle{\frac{\sqrt{k!(k + m)!}}{(k - p)!}e_{k + m}}$ if $ k \geq p+m $$\hfill { } (1.8)$\\

Thus from (1.7), if we denote $\displaystyle{\mathbb{B}_{p} = \{\phi \in \mathbb{B}; \phi(0) = \phi^{'}(0) = ..... = \phi^{(p-1)}(0) = 0 \}}$ then this space is generated by $\displaystyle{\{e_{p}, e_{p+1}, .... \}}$ and the matrix representation of the minimal operator $\mathbb{H}$ generated by the expression $H^{p,m}$ in the basis $\displaystyle{e_{k}; k = p, p+1, ... }$ is given by the symmetric Jacobi matrix $\mathbb{H}$ witch has only two nonzero diagonals. Namely,  its numerical entries are the matrices $\displaystyle{H_{ij}}$ of order $m$ defined by:\\

$\mathbb{J} =$ $ \left\{\begin{array}[c]{l}\displaystyle{H_{i,i} =H_{i,j} = O \quad; if \quad \mid i-j\mid > 1}\quad (i,j = 1, 2, ...)\\
 where\quad O \quad is\quad the \quad zero \quad m\times m\quad matrix\\ \quad \\ and \\\quad\\ \displaystyle{H_{i+1},i =H_{i, i+1}}\quad where\quad H_{i, i+1}\quad is \quad diagonal \quad m\times m\quad matrix\\\quad \\such \quad that\quad its\quad numerical\quad entries \quad are \\\quad \\
\displaystyle{\beta_{k}^{i} = \frac{\sqrt{k!(k + m)!}}{(k - p)!};\quad (i-1)m + 1\leq k \leq im}\quad\\ \end{array}\right. \hfill { }  (1.9)$\\

Let $\displaystyle{A_{i}}$ and $\displaystyle{B_{i} = B_{i}^{*}}$ $(i= 1, 2, .......)$ be $m\times m$ matrices whose entries are complex numbers then the matrix (1.9) is a particular case of the infinite matrix whose general form\\

$\mathfrak{J}=$ $\left(
  \begin{array}{ c c c c c c c c c }
     A_{1} & B_{1} & O & .&.&.&.& \\
     B_{1}^{*}  & A_{2} & B_{2} &\ddots & .&.&.\\
    O &B_{2}^{*} & A_{3} & B_{3}& \ddots &.&.\\
    . & \ddots & \ddots & \ddots & \ddots & \ddots & .&\\
     . & . & \ddots & \ddots & \ddots  & \ddots& O & \\
     . & . & . & \ddots & \ddots  & \ddots & \ddots &\\
     .   &. & . &.& O & \ddots & \ddots& \\
  \end{array} \right).\hfill { }  (1.10)$\\

  \quad\\

  where $O$ is the zero $m\times m$ matrix and the asterisk denotes the adjoint matrix.\\

  Let $l_{m}^{2}(\mathbb{N})$ be the Hilbert space of infinite sequences $\phi = (\phi_{1}, \phi_{2},.....,\phi_{i}, ......)$ with the inner product $\displaystyle{< \phi, \psi > = \sum_{i=1}^{\infty}\phi_{i}\overline{\psi}_{i}}$\\

  where $\displaystyle{\phi_{i} = (\phi_{i}^{1}, \phi_{i}^{2},.....,\phi_{i}^{m}) \in \mathbb{C}^{m}}$ and $\displaystyle{\phi_{i}\overline{\psi}_{i} = \sum_{j=1}^{m}\phi_{i}^{j}\overline{\psi}_{i}^{j}}$\\

  The matrix $\mathfrak{J}$ defines a symmetric operator $\mathfrak{T}$ in $l_{m}^{2}(\mathbb{N})$  according to the formula \\

  $\displaystyle{(\mathfrak{T}\phi)_{i} = B_{i-1}\phi_{i-1} + A_{i}\phi_{i} + B_{i}\phi_{i+1}, i = 1, 2, .....}$ $\hfill { }  (1.11)$\\

where $\displaystyle{\phi_{0} = (\phi_{0}^{1}, \phi_{0}^{2},.....,\phi_{0}^{m}) = (0, 0, .....,0)}$\\

then for our operator, we have:\\

$\displaystyle{(H\phi)_{i} = H_{i,i-1}\phi_{i-1}  + H_{i,i+1}\phi_{i+1}, i = 1, 2, .....}$ $\hfill { }  (1.12)$\\

where $\displaystyle{\phi_{0} = (\phi_{0}^{1}, \phi_{0}^{2},.....,\phi_{0}^{m}) = (0, 0, .....,0)}$\\

\quad\\

{\bf Remark 1.3}\\

1) The closure $\mathbb{T}$ with domain $D(\mathbb{T})$ of the operator $\mathfrak{T}$ is the minimal  closed symmetric operator generated by the expression (1.11) and the boundary condition $\phi_{0} = 0$\\

2) According to Berezanskii, by Chap VII, § 2 [4] , it is well known that the deficiency numbers $n_{+}$ and  $n_{-}$ of the operator $\mathbb{T}$ satisfy the inequalities $0 \leq  n_{+} \leq m$ and $0 \leq  n_{-} \leq m$\\

where  $n_{+}$ is the dimension of $\mathfrak{M}_{z} = (\mathbb{T} - zI)D(\mathbb{T}); \mathfrak{I}m z \neq 0$ and $n_{-}$  is the dimension of the eigensubspace $\mathfrak{N}_{\overline{z}}$ corresponding to the eigenvalue $\overline{z}$ of the operator $\mathbb{T}$\\

3) According to Krein [23] that the operator $\mathbb{T}$ is said completely indeterminate if $n_{+} = n_{-} = m$ and to Kostyuchenko-Mirsoev [22] that the completely indeterminate case holds for the operator $\mathbb{T}$ if and only if all solutions of the vector equation\\

$(\mathfrak{T}\phi)_{i} = z\phi_{i}$ \quad $i=1, 2, .... .$ $\hfill { } (1.13)$\\

for $z  = 0$ belongs to $l_{m}^{2}(\mathbb{N})$\\

In 1949 Krein developed the theory of entire operators with arbitrary finite defect numbers, we refer to [4, 8, 9, 10, 12, 22, 23, 27] and the references therein which are closely connected with this theory. In section 2, we give some properties associated to this theory for the generalized Heun's operator $\displaystyle{H^{p,m} = a^{*^{p}} (a^{m} + a^{*^{m}})a^{p}}$, in particular its completely indeterminacy in Bargmann space.\\
 In section 3, we show that the operators $\displaystyle{\breve{\mathbb{H}} = a^{*^{p}}a^{p+m}\mathbb{U}}$ and $\displaystyle{\breve{\mathbb{H}} +  \breve{\mathbb{H}}^{*}}$;\\ $(p = 1, 2, ....., m=1, 2, .....).$ are chaotic
 where $\displaystyle{\mathbb{U}e_{k} = e_{k+m-1}}$ and
$\displaystyle{\breve{\mathbb{H}}^{*^{p,m}}}$, $\displaystyle{\mathbb{U}^{*}}$ are respectively the adjoint of $\breve{\mathbb{H}}$ and of $\displaystyle{\mathbb{U}}$.{\color{blue}\hf}\\

\begin{center}
\Large{\bf {\color{blue}2. On the completely indeterminacy of generalized Heun operator in Bargmann space}}
\end{center}

In [22], Kostyuchenko and Mirzoev gave some tests for the complete indeterminacy of a Jacobi matrix $\mathfrak{J}$ in terms of entries $A_{i}$ and $B_{i}$ of that matrix. In the following, we give two lemmas witch permit us to show the complete indeterminacy of generalized Heun operator in Bargmann space.\\

For $m = 1, 2, ...$, let $\mathbb{C}^{m}$ be the euclidean $m-$dimentional space and
$\displaystyle{B_{i} = H_{i,i+1}}$ be the diagonal $ m\times m$ matrix such that its numerical entries are \\

$\displaystyle{\beta_{k}^{i} = \frac{\sqrt{k!(k + m)!}}{(k - p)!};\quad (i-1)m + 1\leq k \leq im} (i=1,2, .....)$.\\

By $\mid\mid . \mid\mid$ we denote the spectral matrix norm, then we have:\\

1) $\displaystyle{\mid\mid B_{i} \mid\mid = \frac{\sqrt{(im)![(i+1)m]!}}{(im - p)!} \sim i^{p+\frac{m}{2}}}$ \\

2)  $\displaystyle{\mid\mid B_{i}^{-1} \mid\mid = \frac{1}{\mid\mid B_{i} \mid\mid}}$\\

{\bf Lemma 2.1}\\

Let $\displaystyle{B_{i} = H_{i,i+1}}, (i = 1, 2, ...)$ be the diagonal $ m\times m $ matrix ($m = 1, 2, ...$) such that its numerical entries are \\

$\displaystyle{\beta_{k}^{i} = \frac{\sqrt{k!(k + m)!}}{(k - p)!};\quad (i-1)m + 1\leq k \leq im} (i=1,2, .....)$.\\

then the following inequality holds\\

$\displaystyle{\mid\mid B_{i-1}\mid\mid \mid\mid B_{i+1}\mid\mid \leq \frac{1}{\mid\mid B_{i}^{-1}\mid\mid^{2}}}$\hfill { } (1.16)\\

holds starting from some $i \geq m$\\

{\bf Proof}\\

By using the lemma 1.1 or the behavior of Gamma function $\Gamma(x)$ as \\$\mathcal{R}e x \rightarrow +\infty$ given by Stirling's formula $\displaystyle{\Gamma(x)\sim\sqrt{2\pi}e^{-x}x^{x-\frac{1}{2}}}$, we deduce that \\

$\displaystyle{\beta_{k}^{i} \sim k^{p+\frac{m}{2}}}$ as $k \rightarrow +\infty$.\\

As $\displaystyle{\mid\mid B_{i}\mid\mid = \beta_{im}^{i} = \frac{\sqrt{(im)!(m(i+1))!}}{(im - p)!}\sim (im)^{p+\frac{m}{2}} = (i)^{p+\frac{m}{2}}(m)^{p+\frac{m}{2}}}$

 then  $\displaystyle{\mid\mid B_{i-1}\mid\mid \sim(i-1)^{p+\frac{m}{2}}(m)^{p+\frac{m}{2}}}$, $\displaystyle{\mid\mid B_{i+1}\mid\mid \sim(i+1)^{p+\frac{m}{2}}(m)^{p+\frac{m}{2}}}$ and \\

 $\displaystyle{\mid\mid B_{i-1}\mid\mid}$.$\displaystyle{\mid\mid B_{i+1}\mid\mid}$ $\displaystyle{\sim (i)^{2p+ m}(1- \frac{1}{i^{2}})^{p+\frac{m}{2}}(m)^{2p+m}}$\\

 Now as \\

 $\displaystyle{(1 - \frac{1}{i^{2}})^{p+\frac{m}{2}} \leq 1}$ then \\

 $\displaystyle{\mid\mid B_{i-1}\mid\mid}$ $\displaystyle{\mid\mid B_{i+1}\mid\mid} \leq $ $\displaystyle{\mid\mid B_{i}\mid\mid^{2}}$ and as $\displaystyle{\mid\mid B_{i}^{-1} \mid\mid = \frac{1}{\mid\mid B_{i} \mid\mid}}$ then  (1.16) holds.{\color{blue}\hf}\\

\quad\\

{\bf Lemma 2.2}\\

Let $\displaystyle{B_{i} = H_{i,i+1}}, (i = 1, 2, ...)$ be the diagonal $ m\times m$ matrix , $(m = 1, 2, ...)$ such that its numerical entries are \\

$\displaystyle{\beta_{k}^{i} = \frac{\sqrt{k!(k + m)!}}{(k - p)!};\quad (i-1)m + 1\leq k \leq im}$\\

then if $2p + m > 2$ the following inequality holds\\

$\displaystyle{\sum_{i=1}^{+\infty}\frac{1}{\mid\mid B_{i}\mid\mid} < +\infty}$ $\hfill { } (1.17)$\\

{\bf Proof}\\

As $\displaystyle{\mid\mid B_{i}\mid\mid\sim i^{p + \frac{m}{2}}}$ then  if $p + \frac{m}{2} > 1$ , the serie $\displaystyle{\sum_{i=1}^{+\infty}\frac{1}{i^{p + \frac{m}{2}}}} < +\infty$, it follows that (1.17) holds.{\color{blue}\hf}\\

Now, we prove the following theorem\\

{\bf Theorem 2.3}\\

If $p + \frac{m}{2} > 1$ then the operator $\mathbb{H}$  is completely indeterminate and its deficient numbers satisfy the conditions $n_{+} =  n_{-} = m$.\\

{\bf Proof}\\

By applying the results of Kostyuchenko and Mirsoev [22] to our operator then the completely indeterminate case holds for the opertor $\mathbb{H}$ if and only if all solutions of the vector equation\\

$\displaystyle{ B_{i-1}\phi_{i-1}  + B_{i}\phi_{i+1} = \lambda \phi_{i}\quad (i = 1, 2, .....)}$\\

for $\lambda = 0$ belongs to $l_{m}^{2}(\mathbb{N})$\\

\n where $\displaystyle{B_{i} = H_{i,i+1}}$ is the diagonal $ m\times m$ matrix such that its numerical entries are given by
$\displaystyle{\beta_{k}^{i} = \frac{\sqrt{k!(k + m)!}}{(k - p)!};\quad (i-1)m + 1\leq k \leq im}$\\

Now from (1.12) we consider the system\\

$\displaystyle{ B_{i-1}\phi_{i-1}  + B_{i}\phi_{i+1} = 0 \quad (i = 1, 2, .....)}$\\

where $\displaystyle{\phi_{0} = (\phi_{0}^{1}, \phi_{0}^{2},.....,\phi_{0}^{m}) = (0, 0, .....,0)}$\\

As $\displaystyle{B_{i}^{-1},\quad (i = 1, 2, .....)}$ exist we deduce that the solutions of the above equation have the following explicit form\\

If $i= 2j, \quad (j = 1, 2, ....)$ we have\\

$\displaystyle{\phi_{2j} = 0 }$ and $\displaystyle{\phi_{2j+1} = -(1)^{j}B_{2j}^{-1}B_{2j-1}\times B_{2j-2}^{-1}B_{2j-3} ...... \times B_{2}^{-1}B_{1}\phi_{1}}$\\

This solution belongs to $l_{2}(\mathbb{N})$ if $\displaystyle{\sum_{j=1}^{+\infty}\mid\mid B_{2j}^{-1}B_{2j-1}\times B_{2j-2}^{-1}B_{2j-3} ...... \times B_{2}^{-1}B_{1}\mid\mid^{2} < + \infty}$\\

and\\

if $i= 2j-1, \quad (j = 1, 2, ....)$ we have\\

$\displaystyle{\phi_{2j-1} = 0 }$ and $\displaystyle{\phi_{2j} = -(1)^{j}B_{2j-1}^{-1}B_{2j-2}\times B_{2j-3}^{-1}B_{2j-4} ...... \times B_{3}^{-1}B_{2}\phi_{2}}$\\

This solution belongs to $l_{2}(\mathbb{N})$ if $\displaystyle{\sum_{j=1}^{+\infty}\mid\mid B_{2j-1}^{-1}B_{2j-2}\times B_{2j-3}^{-1}B_{2j-4} ...... \times B_{3}^{-1}B_{2}\mid\mid^{2} < + \infty}$\\

Then the solution generated by the above solutions belongs to $l_{2}(\mathbb{N})$ if \\

$\displaystyle{\sum_{j=1}^{+\infty}\mid\mid B_{2j-1+\epsilon}^{-1}B_{2j-2+\epsilon}\times ...... \times B_{3+\epsilon}^{-1}B_{2+\epsilon}B_{1+\epsilon}^{-1}B_{\epsilon}\mid\mid^{2} < + \infty}$ $\hfill { }    (1.18) $\\

where $\epsilon = 0$ or $\epsilon = 1$ and $\displaystyle{B_{0} = B_{1}^{-1}}$\\

Now as \\

$\displaystyle{\mid\mid B_{2j-1+\epsilon}^{-1}B_{2j-2+\epsilon}\times ...... \times B_{3+\epsilon}^{-1}B_{2+\epsilon}B_{1+\epsilon}^{-1}B_{\epsilon}\mid\mid^{2} \quad \leq}$\\

$\displaystyle{\mid\mid B_{2j-1+\epsilon}^{-1}\mid\mid^{2}\mid\mid B_{2j-2+\epsilon}\mid\mid^{2}\times ...... \times \mid\mid B_{3+\epsilon}^{-1}\mid\mid^{2}\mid\mid B_{2+\epsilon}\mid\mid^{2}\mid\mid B_{1+\epsilon}^{-1}\mid\mid^{2}\mid\mid B_{\epsilon}\mid\mid^{2}}$\\

then it follows from (1.16) of lemma 2.1 that\\

$\displaystyle{\mid\mid B_{2j-1+\epsilon}^{-1}\mid\mid^{2}\times ...... \times \mid\mid B_{3+\epsilon}^{-1}\mid\mid^{2}\mid\mid B_{1+\epsilon}^{-1}\mid\mid^{2} \quad \leq }$\\

$\displaystyle{\frac{1}{\mid\mid B_{2j-2+\epsilon}\mid\mid^{2}\times ...... \times \mid\mid B_{2+\epsilon}\mid\mid^{2}\mid\mid B_{1+\epsilon}^{-1}\mid\mid^{2}\mid\mid B_{\epsilon}\mid\mid \mid\mid B_{2j+\epsilon}\mid\mid }}$\\

and consequently the general term of the series (1.18) do not exceed \\

$\displaystyle{\frac{\mid\mid B_{\epsilon}\mid\mid }{\mid\mid B_{2j+\epsilon}\mid\mid}}$
and (1.17) of the lemma 2.2 ensures the convergence of the series
$\displaystyle{\sum_{j=1}^{+\infty}\frac{1}{\mid\mid B_{2j+\epsilon}\mid\mid}}$.
The  proof of Theorem 2.3 is complete.{\color{blue}\hf}\\

\quad\\

Let $\displaystyle{\breve{\mathbb{H}} = a^{*^{p}}a^{p+m}\mathbb{U}}$ $(p = 1, 2, ....., m=1, 2, .....)$ and $\displaystyle{\breve{\mathbb{H}}^{*}}$ is its adjoint and $\displaystyle{\mathbb{U}e_{k} = e_{k+m-1}}$ . In next section, we study the chaoticity of the operators $\breve{\mathbb{H}}$ and  $\displaystyle{\breve{\mathbb{H}} + \breve{\mathbb{H}}^{*}}$ on Bargmann space in the sense of following Denavey's definition [2], [14] :\\

{\bf Definition 2.4}\\

{\it A linear unbounded densely defined operator $(\mathbb{T}, D(\mathbb{T}))$ on a Banach space $\mathbb{X}$ is called chaotic if the following conditions are met:\\

 1) $\mathbb{T}^{n }$ is closed for all positive integers $n$..\\

 2) there exists an element $\phi \in D(\mathbb{T}^{\infty}) = \cap_{n=1}^{\infty}D(\mathbb{T}^{n})$ whose orbit \\$Orb(\mathbb{T},\phi) = \{\phi, \mathbb{T}\phi, \mathbb{T}^{2}\phi, .....\}$ is dense in $\mathbb{X}$ i.e $\mathbb{T}$ is said to be hyper-cyclic.\\

 3) the set $\{\phi \in \mathbb{X}; \exists \quad j \in\mathbb{ N}$ such that $\mathbb{T}^{j}\phi = \phi\}$ of periodic points of operator $\mathbb{T}$ is dense in $\mathbb{X}$.}\\

{\bf Remark 2.5 }\\

i) It is well known that linear operators in finite-dimensional linear spaces can't
be chaotic but the nonlinear operator may be. Only in infinite-dimensional
linear spaces can linear operators have chaotic properties.\\
These last properties are based on the phenomenon of hypercyclicity or the phenomen of nonwandercity.\\

ii) The study of the phenomenon of hypercyclicity originates in the papers by Birkhoff [7] and Maclane [24] that show, respectively, that the operators of
translation and differentiation, acting on the space of entire functions are hyper-cyclic.\\

iii) Ansari asserts in [1] that powers of a hyper-cyclic bounded operator are also hyper-cyclic\\

iv) For an unbounded operator, Salas exhibit in [25] an unbounded hyper-cyclic
operator whose square is not hyper-cyclic.\\

v) H.N. Salas found in [26] an example of bilateral weighted Shift $\mathbb{T}$ such that both $\mathbb{T}$ and $\mathbb{T}^{*}$ are hypercyclic. The operator $\mathbb{T} \oplus \mathbb{T}^{*}$ is not even cyclic and therefore the direct sum of hypercyclic operators is not always hypercyclic.\\

vi) In Bargmann representation the annihilation operator $a$ is chaotic but $a + a^{*}$ is not chaotic where $a^{*}$ is its adjoint satisfying $[a, a^{*}] = I$.\\

vii) The result of Salas show that one must be careful in the formal manipulation of operators with restricted domains. For such operators it is often more convenient to work with vectors rather than with operators themselves.\\

 \begin{center}
\Large{\bf{\color{blue} 3. On the chaoticity of generalized Heun operator in Bargmann space}}

\end{center}

We begin by recalling some sufficient conditions on hypercyclicity of unbounded operators given by the following (B$\grave{e}$s-Chan-Seubert theorem\\

{\bf Theorem 3.1 } ((B$\grave{e}$s-Chan-Seubert  {\bf [6]}, p.258 )\\

{\it Let $\mathbb{X}$ be  a separable infinite dimensional Banach and let $\mathbb{T}$ be a densely defined linear operator on $\mathbb{X}$. Then $\mathbb{T}$ is hypercyclic if\\
(i) $\mathbb{T}^{m}$ is closed operator for all positive integers $m$.\\
(ii) There exist a dense subset $\mathbb{Y}$ of the domain $D(\mathbb{T})$ of $\mathbb{T}$ and a (possibly nonlinear and discontinuous) mapping $\mathbb{S} : \mathbb{Y} \longrightarrow \mathbb{Y} $ so that $\mathbb{T}\mathbb{S} = I_{\mid \mathbb{Y}}$ ($I_{\mid \mathbb{Y}}$ is identity on $\mathbb{Y}$) and $ \mathbb{T}^{n}, \mathbb{S}^{n} \longrightarrow 0 $ pointwise on $\mathbb{Y}$ as $ n \longrightarrow \infty.$}\\

{\bf Lemma 3.2}\\

Let  $\mathbb{B}$ be the Bargmann space and $a$ and $a^{*}$ are the annihilation and creation operators defined on $\mathbb{B}$ by $a\phi(z) = \phi'(z)$ and $a^{*}\phi(z) = z\phi(z)$ then\\

1) $H^{0,m} = a^{m} + a^{*^{m}}; m = 1, 2$ and $ a^{*^{p}}a^{p}; p = 1, 2, ....$ are not chaotic operators.\\

2) $a$ and $H^{1,1} = a^{*}(a + a^{*})a$ are chaotic operators.\\

{\bf Proof}\\

1) In Bargmann representation we note that:\\

 - The operators $a$ and $a^{*}$ have same domain and that the operators $a^{m} + a^{*^{m}}; m= 1, 2$ are symmetric then they are self-adjoint and consequently they are not chaotic in Bargmann space. We can also use the Carleman'criteria.\\

 - The operators $ a^{*^{p}}a^{p}; p= 1, 2, ....$ are self-adjoint operators with compact resolvent then they are  not chaotic in Bargmann space.\\

2) In [16], it showed that $a^{*^{p}}a^{p+1}$  is chaotic for all $ p \geq 0$ in particular the operator $a$ is chaotic on $\mathbb{B}$ and this result is generalized in [17] to $z^{p}\mathbb{D}^{p+1}$ where $\mathbb{D}$ is Gelfond-Leontiev operator of generalized differentiation [13] acting on generalized Fock-Bargmann space. Recently, for some weighted shift $\mathbb{M}$ defined on $(\Gamma , \chi)$-$Theta$ Fock-Bargmann spaces, we showed in [21] that the operators $\mathbb{M}^{p}\mathbb{D}^{p+1}$ are chaotic for all $ p \geq 0$ where $\mathbb{D}$ is adjoint of $\mathbb{M}$.\\

It is showed in [11] that $H^{1,1} = a^{*}(a + a^{*})a$ is chaotic on $\mathbb{B}$, this last operator play an essential role in Reggeon field theory (see [19] and [20]). Also, we can show that $H^{1,1} = a^{*}(a + a^{*})a$ is chaotic on $\mathbb{B}$ by choosing $\gamma_{n} = \sqrt{n}Logn$  in the theorem [18] below on the chaoticity of the sum of chaotic shifts with their adjoint in Hilbert space.\\

{\bf Theorem 3.3 } [18]\\

Let a linear unbounded densely defined chaotic shift operator $(\mathbb{T}, D(\mathbb{T}))$ on a Hilbert space $\mathbb{E} = \displaystyle{\{\phi; \phi = \sum_{n=1}^{\infty}a_{n}e_{n}\}}$ such that its adjoint is defined by:\\

$\mathbb{T}^{*}e_{n} = \omega_{n}e_{n+1}$ $\hfill { } (3.1)$\\

where $\{e_{n}\}$ is an orthonormal basis of $\mathbb{E}$ and $\omega_{n}$ is positive weight associated to $\mathbb{T}$\\

We assume that\\

(Assumption $Hyp_{1}$)\quad $\displaystyle{\sum_{n=1}^{\infty}\frac{1}{\omega_{n}} < \infty}$ $\hfill { } (3.2)$\\

(Assumption $Hyp_{2}$)\quad $\displaystyle{\omega_{n-1}\omega_{n+1} \leq \omega_{n}^{2}}$ $\hfill { } (3.3)$\\

(Assumption $Hyp_{3})$\quad there exist $\alpha > 0$, $\beta > 0$, $a > 0$, and a sequence $\gamma_{n}$ that:\\

(1) $\displaystyle{\frac{\omega_{n}\gamma_{n}}{\gamma_{n+1}} \geq n^{1+\alpha}}$  $\hfill { } (3.4)$\\

(2) $\displaystyle{\frac{\omega_{n-1}\gamma_{n+1}}{\omega_{n}\gamma_{n-1}} = 1 - \frac{a}{n} + O(\frac{1}{n^{1+\beta}})}$ $\hfill { } (3.5)$\\

and\\

(3) $\displaystyle{\sum_{k=1}^{\infty}\frac{1}{\gamma_{n}^{2}} < \infty}$ $\hfill { } (3.6)$\\

Then  for $\lambda \in \mathbb{C}$ the following recurrence sequence\\

$ (*) \left\{\begin{array}[c]{l}\displaystyle{u_{1}(\lambda) = 1}\\
\displaystyle{u_{2}(\lambda) = \frac{\lambda}{\omega_{1}}} \\
\displaystyle{\omega_{n-1}u_{n-1}(\lambda) + \omega_{n}u_{n+1}(\lambda) = \lambda u_{n}(\lambda)}\\
\end{array}\right.\hfill { }  (3.7)$\\

\quad\\

(i) is solvable  for all $\lambda \in \mathbb{C}$.\\

(ii) $\displaystyle{\sum_{n=1}^{\infty}\mid u_{n}(\lambda)\mid^{2} < \infty}$ for all $\lambda \in \mathbb{C}$.\\

(iii) the spectrum of $\mathbb{T} + \mathbb{T}^{*}$ is the all complex plane $\mathbb{C}$.\\

(iv) $(\mathbb{T} + \mathbb{T}^{*})^{m}$ is closed $\quad \forall \quad m \in \mathbb{N}$.\\

(v) $\mathbb{T} + \mathbb{T}^{*}$ is hypercyclic operator.\\

(vi) $\mathbb{T} + \mathbb{T}^{*}$ is chaotic operator.\\

{\color{blue}\hf}\\

{\bf Remark 3.4}\\

i) In Bargmann representation, the operator $a$ give an example of linear unbounded densely defined chaotic shift operator $(\mathbb{T}, D(\mathbb{T}))$ on a Hilbert space $\mathbb{B}$ such that $\mathbb{T} + \mathbb{T}^{*}$ is not chaotic operator.\\

ii) In Bargmann representation, let $\mathbb{T}$ be linear unbounded densely defined chaotic shift operator such that $n_{+}=n_{-} \neq 0$ where $n_{+}$ and $n_{-}$ are the defect numbers of $\mathbb{T} + \mathbb{T}^{*}$,then the operator $\mathbb{T} + \mathbb{T}^{*}$ is it chaotic ? {\color{blue}\hf}\\

Now we recall the operator $\mathbb{U}$  defined by\\

 $\mathbb{U}e_{k} = e_{k + m -1}$ and its adjoint $\mathbb{U}^{*}e_{k} = e_{k - m +1}; k \geq m$ \\

Let $\displaystyle{\breve{\mathbb{H}} = a^{*^{p}}a^{p+m}\mathbb{U}}$ and $\displaystyle{\breve{\mathbb{H}}^{*} = \mathbb{U}^{*}a^{*^{p+m}}a^{p}}$ then \\

$\displaystyle{\breve{\mathbb{H}}e_{k} = \omega_{k-1}^{p,m}e_{k-1}}$ where \\

$\displaystyle{\omega_{k-1}^{p,m}= \frac{\sqrt{(k-1)!(k-1 + m)!}}{(k-1 -p)!}}$.\\

and\\

$\displaystyle{\breve{\mathbb{H}}^{*}e_{k} = \omega_{k}^{p,m}e_{k+1}}$\\

In the following, we show  that the operator $\displaystyle{\breve{\mathbb{H}}}$ is chaotic and we apply the theorem 3.3 to prove that $\displaystyle{\breve{\mathbb{H}} + \breve{\mathbb{H}}^{*}}$ is also chaotic.\\

{\bf Theorem 3.5}\\

Let $\mathbb{B}_{p}, p = 0, 1, ....$ be the subspace of Bargmann space generated by  $\displaystyle{e_{k}; k \geq p}$ and the operator $\displaystyle{\breve{\mathbb{H}}}$ with domain $\displaystyle{D(\breve{\mathbb{H}}) = \{\phi \in \mathbb{B}_{p}; \breve{\mathbb{H}}\phi \in \mathbb{B}_{p}\}}$ defined by \\

 $\displaystyle{\breve{\mathbb{H}}e_{k} = \omega_{k-1}^{p,m}e_{k-1}}$ where
 $\displaystyle{\omega_{k-1}^{p,m} = \frac{\sqrt{(k-1)!(k-1+m)!}}{(k-p)!}}$\\

 then \\

 $\displaystyle{\breve{\mathbb{H}}}$ is chaotic on $\mathbb{B}_{p}$.\\

 {\bf Proof}\\

To use the theorem of $B\grave{e}s$ and al, we begin by observing that for $\displaystyle{\phi(z) =\sum_{k=p}^{\infty}a_{k}e_{k}(z)}$ such that $\displaystyle{\sum_{k=p}^{\infty}\mid a_{k}\mid^{2} < \infty}$ we have the obvious properties\\

 (i) $\displaystyle{\breve{\mathbb{H}}^{l}\phi(z) =  \sum_{k=p}^{\infty}[\prod_{j=p}^{l+k-1}\omega_{j}^{p,m}]a_{k+l}e_{k}(z)}$ of domain\\

$\displaystyle{D(\breve{\mathbb{H}}^{l}) = \{ \phi = \sum_{k=p}^{\infty}a_{k}e_{k};\sum_{k=p}^{\infty}\mid a_{k}\mid^{2} < \infty \quad and \quad \sum_{k=p}^{\infty}[\prod_{j=p}^{l+k-1}\omega_{j}^{p,m}]^{2}\mid a_{k+l}\mid^{2} < +\infty\}}$\\

 witch is dense in $\displaystyle{\mathbb{B}_{p} \quad\forall\quad l \in \mathbb{N}}$\\

(ii) $\breve{\mathbb{H}}^{l}$ is closed $\quad\forall\quad l \in \mathbb{N}$ and $\breve{\mathbb{H}}^{l} e_{k}(z) = 0$ $\quad\forall\quad l  > k \geq p \geq 0$\\

(iii) As $\omega_{n}^{p,m} \rightarrow +\infty$ then the spectrum of $\breve{\mathbb{H}}$ is all the complex plane.\\

In fact, let $\displaystyle{\phi_{\lambda} = \sum_{k=p}^{\infty}a_{k}e_{k}}$ with $\displaystyle{a_{k} = \prod_{j=p}^{k-1}\frac{\lambda}{\omega_{j}^{p,m}}}$ i.e \\

$\displaystyle{\phi_{\lambda} =\sum_{k=p}^{+\infty}[\prod_{j=p}^{k-1}\frac{\lambda}{\omega_{j}^{p,m}}]e_{k}}$ then as $a_{p} = 0$ we deduce that\\

$\breve{\mathbb{H}}\phi_{\lambda} = \lambda\phi_{\lambda} , \forall \quad \lambda \in \mathbb{C}$.$\hfill { }  (3.8)$\\\\

and as
$\displaystyle{\sum_{k=p}^{\infty}[\prod_{j=p}^{k}\frac{\lambda}{\omega_{j}^{p,m}}]^{2} < +\infty}$
then $\phi_{\lambda} \in D(\breve{\mathbb{H}})$\\

Now, take $\mathbb{Y}$ the linear subspace generated by finite combinations of basis $\{e_{k}\}_{k=p}^{\infty}$, this subspace $\mathbb{Y}$ is dense in $\displaystyle{\mathbb{B}_{p}}$ and we define on it the operator $\mathbb{S}$ acting on $\displaystyle{\phi = \sum_{k=p}^{N}a_{k}e_{k}}$ as following \\

$\displaystyle{\mathbb{S}\phi = \sum_{k=p}^{N+1}\frac{a_{k-1}}{\omega_{k-1}^{p,m}}e_{k}}$ $\hfill { }  (4.5)$\\ then \\

$\displaystyle{\mathbb{S}^{n}e_{k} = \frac{1}{\prod_{j=k}^{n+k}\omega_{j}^{p,m}}e_{k+n}}$\\

as $\displaystyle{\prod_{j=p}^{n}\omega_{j}^{p,m} \rightarrow +\infty}$ as $n \rightarrow +\infty$ we get\\

$\displaystyle{\mathbb{S}^{n}e_{k} \rightarrow 0}$ in $\displaystyle{\mathbb{B}_{p}}$ as $n \rightarrow +\infty$ $\hfill { }  (3.9)$\\

By noting that $\breve{\mathbb{H}}^{n}e_{k} = 0 $ for $n > k$ and any element of $\mathbb{Y}$ can be annihilated by a finite power of $\breve{\mathbb{H}}$ and $\breve{\mathbb{H}}\mathbb{S} = \mathbb{I}_{\mid \mathbb{Y}}$ then the hyperciclycity of $\breve{\mathbb{H}}$ follows from the theorem of $ B\grave{e}s$ and al. recalled above.{\color{blue}\hf}\\

We shall now show that $\breve{\mathbb{H}}$ has a dense set of periodic points.\\

To see this, it suffices to show that for every element $\phi$ in the dense subspace $\mathbb{Y}$ there is a periodic point $\psi$ arbitrarily close to it.\\

For $s \geq p$ and $N \geq s$ we put\\

$\displaystyle{\varphi_{s,N}(z) = e_{s}(z) + \sum_{k=s+1}^{\infty}[\prod_{j=s}^{kN + s -1}\frac{1}{\omega_{j}^{p,m}}]e_{kN +s}(z)}$ $\hfill { }  (3.10)$\\

Then we have the following obvious lemma\\

{\bf Lemma 3.6}\\

(i) $\displaystyle{\breve{\mathbb{H}}^{N}\prod_{j=0}^{kN -1}\frac{1}{\omega_{j}^{p,m}} e_{kN} =\prod_{j=0}^{(k-1)N -1}\frac{1}{\omega_{j}^{p,m}} e_{(k-1)N} \quad \forall \quad k \geq p}$\\

(ii) $\displaystyle{\breve{\mathbb{H}}^{N}\prod_{j=s}^{kN -1+s}\frac{1}{\omega_{j}^{p,m}}e_{kN+s}} = \displaystyle{\prod_{j=s}^{(k-1)N -1}\frac{1}{\omega_{j}^{p,m}} e_{(k-1)N+s}}$ for $s\geq p$,$N\geq s$ and $k \geq p$\\

(iii) $\displaystyle{\varphi_{s,N}}$ is $N$-periodic point of  $\breve{\mathbb{H}}$.\\

(iv) $\displaystyle{\varphi_{s,N} \in D(\breve{\mathbb{H}}^{N})}$.{\color{blue}\hf}\\

Now, Let \\

$\phi(z) = \displaystyle{\sum_{s=p}^{M}a_{s}e_{s}(z)}$ $\hfill { }  (3.11)$\\

such that \\

$\displaystyle{\mid a_{s}\prod_{j=p}^{s-1}\omega_{j}^{p,m}\mid < 1 ; s= p, p+1, ........, M}$ $\hfill { }  (3.12)$\\

and we choose the periodic point for $\breve{\mathbb{H}}$ as $\psi(z)$\\

 $\psi(z) = \displaystyle{\sum_{s=p}^{M}a_{s}\varphi_{s,N}(z)}$ $\hfill { }  (3.13)$\\

then there exists an $N \geq M$ such that\\

$\mid\mid \phi - \psi \mid\mid \leq \epsilon \quad \forall \quad \epsilon > 0$. $\hfill { }  (3.14)$\\

\quad\\

{\bf Remark 3.7}\\

We can also use the results of Bermudez et al [5] to prove the chaoticity of our operator $\displaystyle{\breve{\mathbb{H}}}$ is chotic.\\

{\bf Theorem 3.8}\\

Let $\mathbb{B}_{p}$ the subspace of Bargmann space generated by  $\displaystyle{e_{k}; k \geq p}$ then \\

$\displaystyle{\breve{\mathbb{H}}} + \displaystyle{\breve{\mathbb{H}}^{*}}$ is chaotic
where $\displaystyle{\breve{\mathbb{H}}^{*}}$ is adjoint operator of $\displaystyle{\breve{\mathbb{H}}}$\\

{\bf Proof }\\

For $\displaystyle{\omega_{k}^{p,m} = \sqrt{k!}\frac{\sqrt{(k+m)!}}{(k-m)!} \sim k^{p+\frac{m}{2}}}$ the assumptions (3.2), (3.3) and (3.6) of theorem 3.3 hold. It remains to check the validity of assumptions (3.4) and (3.5)\\

i) The assumption (3.4) is valid. In fact, as $\displaystyle{\omega_{k}^{p,m} = \sqrt{k!}\frac{\sqrt{(k+m)!}}{(k-m)!} \sim k^{p+\frac{m}{2}}}$ then \\

$\displaystyle{\frac{\omega_{k-1}^{p,m}}{\omega_{k}^{p,m}} = (1 - \frac{1}{k})^{p+\frac{m}{2}} = 1 - \frac{p+\frac{m}{2}}{k} + O(\frac{1}{k^{2}})}$\\

ii)The assumption (3.5) is valid. In fact, if we choose $\displaystyle{\gamma_{k} = k^{\frac{m}{2}}Log(k^{p}) = pk^{\frac{m}{2}}Log(k)}$ then\\

$\displaystyle{\frac{\gamma_{k}}{\gamma_{k+1}} = 1 - \frac{1 + \frac{m}{2}}{k} + O(\frac{1}{k^{2}})}$\\

and\\

for $p + \frac{m}{2} > 2$ we get\\

$\displaystyle{\omega_{k}\frac{\gamma_{k}}{\gamma_{k+1}} \geq k^{1+\beta}}$ where $1 +\beta = p + \frac{m}{2} - 1$.{\color{blue}\hf}\\

\newpage
\begin{center}
\Large{\bf{\color{blue} References}}
\end{center}

\n{\bf [1] } S.I. Ansari, Hypercyclic and cyclic vectors. J. Funct. Anal. 128(2),(1995), 374-383.\\

\n{\bf [2] } J. Banks, J. Brooks, G.Cairns, G. Davis and P.Stacey, On Devaney's definition of chaos, Amer. Math. Monthly 99, (1992), 332-334.\\

\n{\bf [3] } V. Bargmann, On Hilbert space of analytic functions and associated integral transform, Part I, Commun. Pure App. Math., 14,(1961), 187-214.\\

\n{\bf [4] } Yu.M. Berezanskii, Expansion in eigenfunctions of selfadjoint operators. Providence, RI: Am. Math.Soc., (1968)\\

\n{\bf [5] } T. Bermudez, A. Bonilla and J. L. Torrea, chaotic behavior of the Riesz transforms for Hermite expansions, J. Math. Anal. Appl., 337, (2008), 702-711.\\

\n{\bf [6] } J. B$\grave{e}$s, K. Chan and S. Seubert, Chaotic unbounded differentiation operators, Integral Equations Operators Theory, 40, (2001), 257-267.\\

\n{\bf [7] } G. D. Birkhoff, D$\acute{e}$monstration  d'un th$\acute{e}$or$\grave{e}$me $ \acute{e}$l$\acute{e}$mentaire sur les fonctions  enti$\grave{e}$res, C. R. Acad. Sci. Paris 189, (1929), 473-475.\\

\n{\bf [8] } A. L. Chistyakov, Deficiency numbers of symmetric operators in a direct sum of Hilbert spaces," I, Vestn. Mosk. Univ., Ser. 1, Mat. Mekh., No. 3, 5-21 (1969); II, No. 4, 3-5 (1969) 3-5.\\

\n{\bf [9] } A. L. Chistyakov, "Deficiency indices of $J_{m}$-matrices and differential operators with polynomial coefficients, Mat. Sb., No. 4,  (1971) 474-503.\\

\n{\bf [10] }  M. L. Gorbachuk, and  V. I. Gorbachuk, Krein's lectures on entire operators. Operator Theory: Advances and Applications, 97. Birkha$\ddot{u}$ser Verlag, Basel, (1997).\\

\n{\bf [11] } A. Decarreau, H. Emamirad, A. Intissar, Chaoticit$\acute{e}$ de l'op$\acute{e}$rateur de Gribov dans l'espace de Bargmann, C. R. Acad. Sci. Paris , 331, (2000)\\

\n{\bf [12] } A. Devinatz, The deficciency index problem for ordinary selfadjoint differential operators, Bulletin of the American Mathematical Society, Volume 79, Number 6, November (1973), 1109-1127\\

\n{\bf [13] } A.O. Gelfond and A.F.Leontiev , Mat. Sb. 29 (3), (1951), 477-500.\\

\n{\bf [14] } A. Gulisashvili and C. MacCluer, linear chaos in the unforced quantum harmonic oscillator, J. Dyn. Syst. Meeas, Control, 118, (1996),337-338.\\

\n{\bf [15] } A.  Intissar, Analyse de scattering d'un op$\acute{e}$rateur cubique de Heun dans l'espace de Bargmann, Comm. Math.Phys. 199, (1998), 243-256.\\

\n{\bf [16] } A. Intissar, On a chaotic weighted shift $z^{p}\frac{d^{p+1}}{dz^{p+1}}$ of
order $p$ in Bargmann space, Advances in Mathematical Physics, (2011), Article ID 471314, 11 pages.\\

\n{\bf [17] } A. Intissar, On a chaotic weighted shift $z^{p}\mathbb{D}^{p+1}$
 of order $p$ in generalized Fock-Bargmann spaces, Mathematica Aeterna, V. 3, no 7, (2013), 519-534.\\

\n{\bf [18] } A. Intissar, On chaoticity of the sum of chaotic shifts with their adjoints
in Hilbert space and applications to some chaotic weighted shifts acting on some Fock-Bargmann spaces, arXiv. 1311.1394v1, submitted to Advances in Theoretical and Mathematical Physics, (2013)\\

\n{\bf [19] } A. Intissar, Etude spectrale d'une famille d'op\'erateurs non sym\'etriques intervenant dans la th\'eorie des champs des reggeons, Comm. Math.Phys. 113 (2), (1987), 263-297\\

\n{\bf [20] } A. Intissar, Spectral analysis of non-selfadjoint Jacobi-Gribov operator
and asymptotic analysis of its generalized eigenvectors, accepted for publication in Advances in Mathematics (China), (2014).\\

\n{\bf [21] } A. Intissar, A short note on the chaoticity of a weight shift on concrete orthonormal basis associated to some Fock-Bargmann space, Journal of Mathematical Physics, 55, 011502 (2014)\\

\n{\bf [22] } A.G.Kostyuchenko and K.A. Mirsoev: Three-term recurrence relations with matrix coefficients, the completely indeterminate case, Mat.Zametki, 63, no. 5 (1998) 709-716 \\

\n{\bf [23] } M.G. Krein, Infinite J-matrices and the matrix moment problem, Dokl. Akad. Nauk SSSR, 69, no. 3 (1949) 125-128)\\

\n{\bf [24] } G. R. Maclane, Sequences of derivatives and normal families, J. Anal.
Math. 2, (1952)\\

\n{\bf [25] } H.N. Salas, Pathological hypercyclic operators, Arch. Math. 86, (2006), 241-250.\\

\n{\bf [26] } H.N. Salas, A hypercyclic operator whose adjoint is also hypercyclic, Proc. Amer. Math. Soc, ********.\\

\n{\bf [27] } L.O. Silva and J. H. Toloza, The class of n-entire operators, J. Phys. A 46,
025202 (23 pp), (2013).\\

\end{document}